\documentclass{article}

\usepackage{graphicx,amsmath,amsfonts,amssymb,url,booktabs,cite,authblk,epstopdf,multirow}
\usepackage[hang,flushmargin]{footmisc} 
\usepackage{graphicx}
\usepackage{amsmath}
\usepackage{amsthm}
\usepackage{amsfonts}
\usepackage{amssymb}
\usepackage{subfig}
\usepackage{url}
\usepackage{verbatim}
\usepackage{color}
\definecolor{red}{cmyk}{0,1,1,0}

\definecolor{blue}{cmyk}{1,0,0,0}
\def\blue{\color{blue}}
\definecolor{gray}{rgb}{0.7,0.7,0.7}

\definecolor{green}{rgb}{0,1,0}

\usepackage[mathlines]{lineno}

\newtheorem{thm}{Theorem}[section]
\newtheorem{defn}[thm]{Definition}
\newtheorem{prop}[thm]{Proposition}
\newtheorem{cor}[thm]{Corollary}
\newtheorem{lem}[thm]{Lemma}

\newtheorem{ex}[thm]{Example}
\newtheorem{quest}[thm]{Question}
\newtheorem{obs}[thm]{Observation}

\newtheorem{rem}[thm]{Remark}

\def\mtx#1{\begin{bmatrix} #1 \end{bmatrix}}
\def\ord#1{| #1 |} 

\newcommand{\R}{\mathbb{R}}

\newcommand{\clf}{\mathcal{F}}

\newcommand{\Rnn}{\R^{n\times n}}

\newcommand{\SG}{\mathcal{S}(G)}
\newcommand{\sym}{\mathcal{S}}
\newcommand{\ZT}{\operatorname{Eff}}
\newcommand{\FT}{\mathcal{F}_{eff}}

\newcommand{\nul}{\operatorname{null}}

\newcommand{\M}{\operatorname{M}}

\newcommand{\ZFN}{\operatorname{Z}}
\newcommand{\PC}{\operatorname{P}}

\newcommand{\T}{\operatorname{Term}}
\newcommand{\Rev}{\operatorname{Rev}}

\newcommand{\x}{\times}

\newcommand{\diam}{\operatorname{diam}}
\newcommand{\pt}{\operatorname{pt}}
\newcommand{\len}{\operatorname{len}}

\newcommand{\PT}{\operatorname{PT}}
\newcommand{\pd}{\operatorname{pd}}
\newcommand{\prev}{\operatorname{prev}}
\newcommand{\next}{\operatorname{next}}
\newcommand{\first}{\operatorname{first}}
\newcommand{\last}{\operatorname{last}}
\newcommand{\pth}{\operatorname{path}}
\newcommand{\npath}{\overline{\pth}}

\newcommand{\bit}{\begin{itemize}}
\newcommand{\eit}{\end{itemize}}
\newcommand{\ben}{\begin{enumerate}}
\newcommand{\een}{\end{enumerate}}
\newcommand{\beq}{\begin{equation}}
\newcommand{\eeq}{\end{equation}}
\newcommand{\bea}{\begin{eqnarray*}}
\newcommand{\eea}{\end{eqnarray*}}
\newcommand{\bpf}{\begin{proof}}
\newcommand{\epf}{\end{proof}\ms}
\newcommand{\bal}{\begin{align*}}
\newcommand{\eal}{\end{align*}}
\newcommand{\ms}{\medskip}
\newcommand{\noi}{\noindent}
\newcommand{\dep}{\displaystyle}
\newcommand{\bmt}{\begin{pmatrix}}
\newcommand{\emt}{\end{pmatrix}}

\title{Propagation time for zero forcing on a graph}
\author{Leslie Hogben\thanks{Department of Mathematics, Iowa State University, Ames, IA 50011, USA (lhogben@iastate.edu) and American Institute of Mathematics, 360 Portage Ave, Palo Alto, CA 94306, USA (hogben@aimath.org).}, My Huynh\thanks{Department of Mathematics, Arizona State University, Tempe, AZ 85287 (mthuynh1@asu.edu).  Research  supported by DMS 0502354 and DMS 0750986.}, Nicole Kingsley\thanks{Department of Mathematics, Iowa State University, Ames, IA 50011, USA (nkingsle@iastate.edu).}, Sarah Meyer\thanks{Department of Mathematics, Smith College, Northampton, MA 01063, USA (smeyer@smith.edu). Research  supported by DMS 0750986.}, Shanise Walker\thanks{Department of Mathematics, University of Georgia, Athens, GA 30602, USA (shanise1@uga.edu).  Research  supported by DMS 0502354 and DMS 0750986}, Michael Young\thanks{Department of Mathematics, Iowa State University, Ames, IA 50011, USA (myoung@iastate.edu).  Research  supported by DMS 0946431.}}

\begin{document}
\maketitle

\begin{abstract}   Zero forcing (also called graph infection) on a simple, undirected graph $G$ is based on the color-change rule: If each vertex of $G$ is colored either white or black, and vertex $v$ is a black vertex with only one white neighbor $w$, then change the color of $w$ to  black. A minimum zero forcing set is a set of black vertices of minimum cardinality
that can color the entire graph black using the color change rule. The propagation time of a zero forcing set $B$ of graph $G$ is the minimum number of steps that it takes to force all the vertices of $G$ black, starting with the vertices in $B$ black and  performing independent forces simultaneously. The minimum and maximum propagation times of a graph are taken over all minimum zero forcing sets of the graph.
 It is shown that a connected graph of order at least two has more than one minimum zero forcing set realizing minimum propagation time.  Graphs $G$ having extreme minimum  propagation times $|G| - 1$, $|G| - 2$, and $0$ are characterized, and results
regarding graphs having minimum propagation time $1$ are established.  It is shown that the diameter is an upper bound for maximum propagation time for a tree, but in general  propagation time and diameter
of a graph are not comparable. 
 \end{abstract}

\noi{\bf Keywords}
 zero forcing number, propagation time, graph \ms

\noi{\bf AMS subject classification}
 05C50, 05C12, 05C15, 05C57, 81Q93, 82B20, 82C20

\section{Propagation time}

All graphs are simple, finite, and undirected.
In a graph $G$ where some vertices are colored black and the remaining vertices are colored white, the {\em color change rule} is:
If $v$ is black and $w$ is the only white neighbor of  $v$, then change the color of $w$ to black; if we apply the color change rule to $v$ to change the color of $w$, we say $v$ {\em forces} $w$ and write  $v\to w$ (note that there may be a choice involved, since as we record forces, only one vertex actually forces $w$, but more than one may be able to). Given an initial set $B$ of black vertices, the {\em final coloring} of $B$ is the set of black vertices that results from applying the  color change rule until no more changes are possible.  For a given graph $G$ and set  of vertices $B$, the final coloring is unique \cite{AIM}.  A {\em  zero forcing set} is an initial set $B$ of vertices such that the final coloring of $B$ is  $V(G)$.  A {\em minimum zero forcing set} of a  graph $G$   is a  zero forcing set  of $G$ of   minimum cardinality, and  the  {\em   zero forcing number}, denoted $\ZFN(G)$,  is the cardinality of a minimum zero forcing set.

Zero forcing, also known as graph infection or graph propagation, was introduced independently in \cite{AIM} for  study of minimum rank problems in combinatorial matrix theory,  and in \cite{graphinfect} for study of control of quantum systems.  Propagation time of a zero forcing set, which describes the number of steps needed to fully color a graph performing independent forces simultaneously, was implicit in  \cite{graphinfect} and explicit in  \cite{Sev}.  Recently Chilakamarri et al. \cite{CDKY} determined the  propagation time, which they call the iteration index, for a  number of families of graphs including Cartesian products and various grid graphs.
Control of an entire network by sequential operations on a subset of particles is valuable \cite{Sev} and the number of steps needed to obtain this control (propagation time) is a significant part of the process.  In this paper we systematically study propagation time.  
\begin{defn} \label{ztdefn}{\rm Let $G=(V,E)$ be a graph and $B$ a zero forcing set of $G$. Define $B^{(0)}=B$, and for $t \ge 0$,  $B^{(t+1)}$ is the set of vertices  $w$ for which there exists a vertex $ b \in \bigcup_{s=0}^{t} B^{(s)}$ such that $w$ is the only neighbor of $b$ not in $\bigcup_{s=0}^{t} B^{(s)}$. The {\em propagation  time of $B$} in $G$, denoted  $\pt(G,B)$,  is the smallest integer $t_0$ such that $V=\bigcup_{t=0}^{t_0} B^{(t)}$. 
}
\end{defn}

Two minimum zero forcing sets of the same graph may have different propagation times, as the following example illustrates.

  \begin{figure}[!ht] 
\begin{center}\scalebox{.35}{\includegraphics{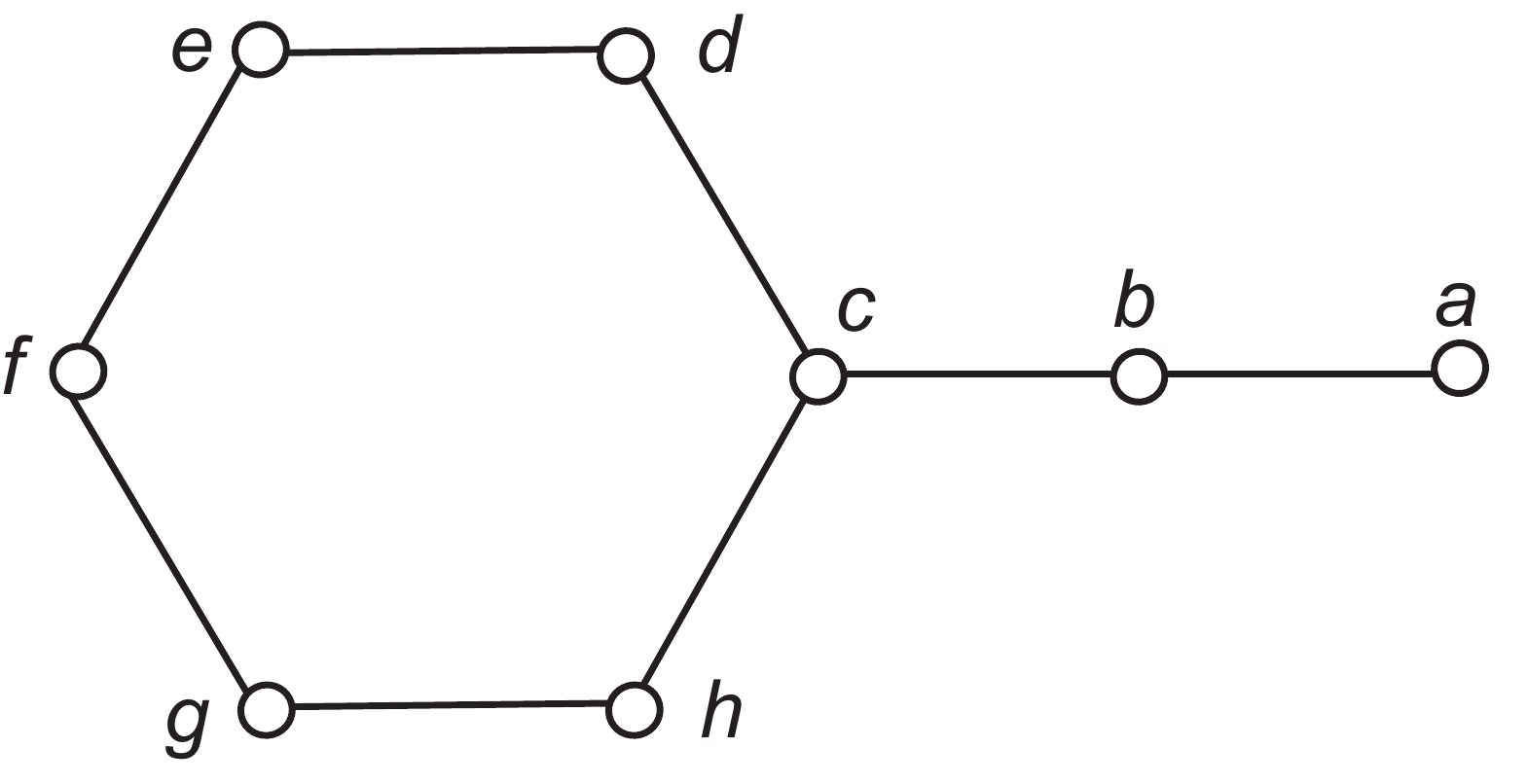}}
\caption{The graph $G$ for Example \ref{ex1} 
} \label{reverseex}
\end{center}
\end{figure}

\begin{ex} \label{ex1}{\rm Let $G$ be the graph in Figure \ref{reverseex}. Let $B_1=\{g,h\}$ and $B_2=\{a,d\}$. Then ${B_1}^{(1)}=\{f,c\}$, ${B_1}^{(2)}=\{e\}$, ${B_1}^{(3)}=\{d\}$, ${B_1}^{(4)}=\{b\}$, and ${B_1}^{(5)}=\{a\}$, so $\pt(G,B_1)=5$.  However, ${B_2}^{(1)}=\{b\}$, ${B_2}^{(2)}=\{c\}$, ${B_2}^{(3)}=\{e,h\}$, and ${B_2}^{(4)}=\{f,g\}$, so $\pt(G,B_2)=4$.}
\end{ex}

\begin{defn}{\rm The {\em minimum propagation  time} of $G$ is 
\[\pt(G)=\min\{\pt(G,B)\,|\,B\text{ is a minimum zero forcing set of } G\}. \] }
\end{defn}

\begin{defn}{\rm Two minimum zero forcing sets $B_1$ and $B_2$ of a graph $G$ are {\em isomorphic} if there is a graph automorphism $\varphi$ of $G$ such that $\varphi(B_1)=B_2$.}
\end{defn}

It is obvious that isomorphic zero forcing sets have the same propagation time, but  a graph  may have  non-isomorphic minimum zero forcing sets and  have the property that all minimum zero forcing sets have the same propagation time. 

\begin{ex}\label{dartex}{\rm Up to isomorphism, the minimum zero forcing sets of the dart shown in Figure \ref{dart} are $\{a,c\}$, $\{b,c\},$ and $\{c,d\}$. Each of these sets has propagation time $3$. 
 \begin{figure}[!ht] 
\begin{center}\scalebox{.35}{\includegraphics{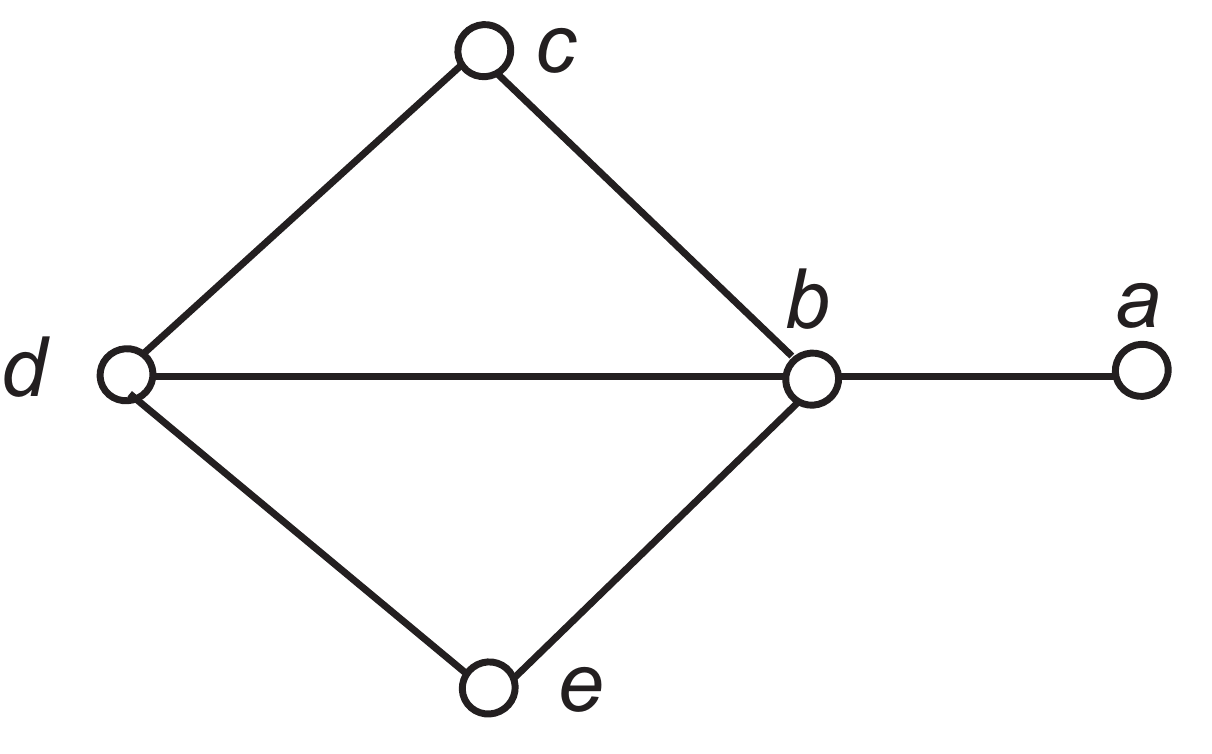}}
\caption{The dart} \label{dart}
\end{center}
\end{figure}
}
\end{ex}

 The   minimum propagation time of a graph $G$ is not subgraph monotone. For example, it is easy to see  that the 4-cycle has $\ZFN(C_4) = 2$ and $\pt(C_4) = 1$. By deleting one edge of $C_4$, it  becomes a path $P_3$, which has $\ZFN(P_3) = 1$ and $\pt(P_3) = 2$. 

A minimum zero forcing set that achieves minimum propagation time plays a central role in our study, and we name such a set.

 \begin{defn} {\rm 
 A    subset $B$ of vertices of $G$ is an {\em efficient zero forcing set} for $G$ if $B$ is a  minimum zero forcing set of $G$ and $\pt(G,B)=\pt(G)$. Define \[\ZT(G) = \{ B \, | \, B \text{ is an efficient zero forcing set of } G\}.\]
}\end{defn}

We can also consider maximum propagation time.

\begin{defn} \label{ztmaxdefn} {\rm The {\em maximum propagation  time} of $G$ is defined as \[\PT(G)=\max\{\pt(G,B)\,|\,B\text{ is a minimum zero forcing set of } G\}. \] }
\end{defn}

The bounds in the next remark were also observed in \cite{CDKY}.

\begin{rem}\label{ptPTbds} {\rm Let $G$ be a graph.  Then
\[\PT(G)  \le  \ord G - \ZFN(G)\] because at least one force must be performed at each time step, and
\[ \frac {\ord G-\ZFN(G)}{\ZFN(G)}   \le  \pt(G)\] because using a given zero forcing set $B$, at most $\ord B$ forces can be performed at any one time step.
}\end{rem}

\begin{defn} \label{intervaldefn} {\rm The {\em propagation time interval} of $G$ is defined as \[[\pt(G),\PT(G)]=\{\pt(G),\pt(G)+1,\dots,\PT(G)-1,\PT(G)\}. \]  The {\em propagation time discrepancy} of $G$ is defined as \[\pd(G) = \PT(G) - \pt(G).\]}
\end{defn}

\begin{figure}[!ht] 
\begin{center}\scalebox{.4}{\includegraphics{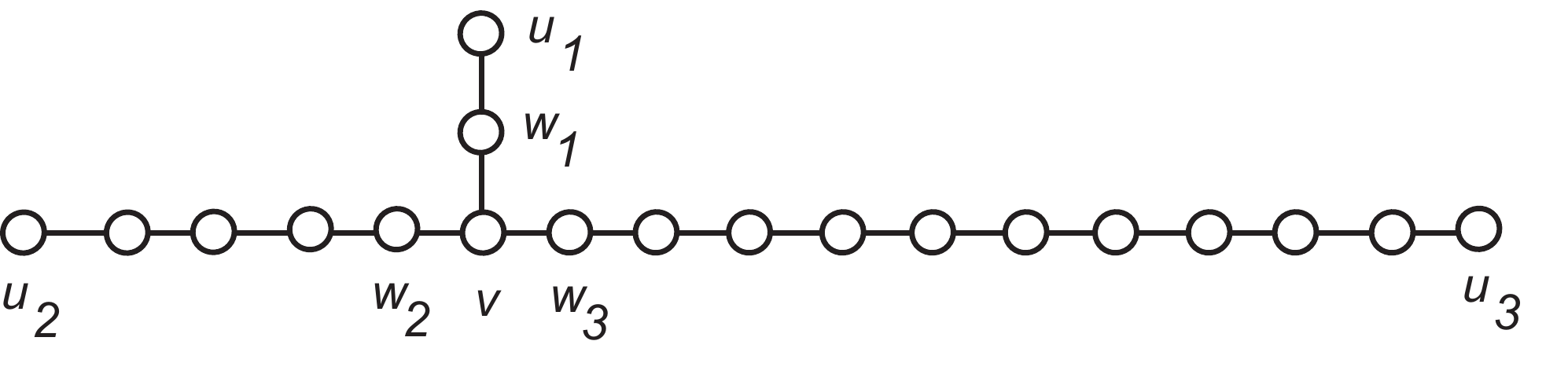}}
\caption{The tree $S(2,5,11)$} \label{s2511}
\end{center}\vspace{-5mm}
\end{figure}

It is not the case that every integer in the propagation time interval is the propagation time of a minimum zero forcing set; this can be seen in the next example. 
Let $S(e_1,e_2,e_3)$ be the generalized star with three arms having $e_1,e_2,e_3$ vertices with $e_1\le e_2\le e_3$; $S(2,5,11)$ is shown  in Figure \ref{s2511}.  

\begin{ex}\label{exT} {\rm  Consider $S(e_1,e_2,e_3)$ with  $1<e_1<e_2<e_3$.  The vertices of degree one are denoted by $u_1, u_2, u_3$, the vertex of degree three is denoted by $v$, and neighbors of $v$ are denoted by $w_1, w_2, w_3$.  The minimum zero forcing sets and their propagation times are shown in Table \ref{tabS2511}.   Observe that  the propagation time interval of $S(2,5,11)$ is $[12,16]$,  but there is no minimum zero forcing set with propagation time $14$. The propagation discrepancy is $\pd(S(2,5,11))=4$.

    \begin{table} [h!]
        \caption{Minimum zero forcing sets and propagation times of $S(e_1,e_2,e_3)$}\label{tabS2511}
        \centering
  \begin{tabular}{|c|c|c|}
        \hline
       $B$   & $\pt(S(e_1,e_2,e_3),B)$ & $\pt(S(2,5,11),B)$ \\
      \hline\hline
  $\{u_1,u_2\}$   & $e_2+e_3-1$ & $15$ \\
        \hline
          $\{u_3,w_2\}$ & $e_2+e_3-1$ & $15$  \\
        \hline
            $\{u_3, w_1\}$  & $e_2+e_3$ & $16$ \\
        \hline
               $\{u_1,u_3\}$  & $e_2+e_3-1$ & $15$ \\
        \hline
          $\{u_2,w_3\}$  & $e_2+e_3-1$ & $15$ \\
        \hline
                 $\{u_2,w_1\}$  & $e_2+e_3$ & $16$ \\
           \hline
                $\{u_2,u_3\}$  & $e_1+e_3-1$ & $12$ \\
       \hline
                 $\{u_1,w_3\}$ & $e_1+e_3-1$ & $12$ \\
        \hline
                    $\{u_1,w_2\}$  & $e_1+e_3$ & $13$ \\
        \hline
          \end{tabular}
      \end{table}
 }\end{ex}

The next remark provides a necessary condition for a graph $G$ to have $\pd(G)=0$.
\begin{rem}{\rm Let $G$ be a graph. Then every minimum zero forcing set of $G$ is an efficient zero forcing set if and only if $\pd(G) = 0$.  In \cite{Zplus}, it is proven that the intersection of all minimum zero forcing sets is the empty set. Hence, $\pd(G)=0$ implies $\dep \bigcap_{B \in \ZT(G)} B = \emptyset$. 
}
\end{rem}

 In Section \ref{sunique} we establish properties of efficient zero forcing sets, including that  no connected graph of order more than one has a unique efficient zero forcing set.  
  In Section \ref{sextreme} we characterize graphs having extreme  propagation time.  In Section \ref{sdiam} we examine the relationship between propagation time and diameter.


\section{Efficient zero forcing sets} 
\label{sunique}

In \cite{Zplus} it was shown that for a connected graph of order at least two, there must be more than one minimum zero forcing set and furthermore, no vertex is in every minimum zero forcing set.  This raises the questions of whether the analogous properties are true for efficient zero forcing sets (Questions \ref{unique} and \ref{intersect} below).

\begin{quest}\label{unique} Is there a connected graph of order at least two that has a unique efficient zero forcing set? \end{quest} 

We  show that the answer to Question \ref{unique} is negative.  First we need some terminology.  For a given zero forcing set $B$ of $G$, construct the  final coloring, listing the forces in the order in which they were performed.  This list is a {\em chronological list of forces} of $B$ \cite{cancun}.  
Many definitions and results concerning lists of forces that have appeared in the literature involve  chronological (ordered) lists of forces.   For the study of propagation time, the order of forces is often dictated by performing a force as soon as possible (propagating). Thus unordered sets of forces are more useful than ordered lists when studying propagation time, and we extend terminology from chronological lists of forces to sets of forces.  

\begin{defn} {\rm  Let $G=(V,E)$ be a graph, $B$ a zero forcing set of $G$.  The unordered set of forces in a chronological list of forces of $B$ is called {\em a set of forces} of $B$.   
}
\end{defn}

Observe that if $B$ is a zero forcing set and $\clf$ is a set of forces of $B$, then the cardinality of $\clf$ is $\ord G-\ord B$.  The ideas of terminus and reverse set of forces, introduced  in \cite{Zplus} for a chronological list of forces and defined below for a set of forces, are used to answer Question \ref{unique} negatively (by constructing the  terminus of a set of forces of an efficient zero forcing set).   

\begin{defn} {\rm  Let $G$ be a graph, let $B$ be a zero forcing set of $G$, and let $\clf$ be a set of forces of $B$.
The {\em terminus} of   $\clf$, denoted  $\T(\clf)$, is the set of vertices that do not perform a force in $\clf$. 
The {\em reverse set of forces} of $\clf$, denoted here as $\Rev(\clf)$, is the result of reversing each force in $\clf$.
A  {\em forcing chain} of  $\clf$ is a sequence of vertices $(v_1,v_2,\dots,v_k)$ such that for $i=1,\dots,k-1$,  $v_i $ forces $v_{i+1}$ in $\clf$ ($k=1$ is permitted).    A {\em maximal forcing chain} is a forcing chain  that is not a proper subset of another  forcing chain. 
}\end{defn}

The name ``terminus" reflects the fact that a vertex does not perform a force in $\clf$ if and only if it is the end point of a maximal forcing chain (the latter is the definition used in \cite{Zplus}, where such a set is called a reversal of $B$). 
  In \cite{Zplus}, it is shown that if $B$ is a zero forcing set of $G$ and $\clf$ is a chronological list of forces, then the terminus of $\clf$ is also a zero forcing set of $G$, with the reverse chronological list of forces (to construct a reverse chronological list of forces of $\clf$,  
  write the chronological list of forces in reverse order and reverse each force in $\clf$). 
   \begin{obs}\label{obsterm} Let $G$ be a graph, $B$ a minimum zero forcing set of $G$, and $\clf$ a set of forces of $B$. Then $\Rev(\clf)$ is a  set of forces  of $\T(\clf)$ and $B=\T(\Rev(\clf))$.
\end{obs}
 
When studying propagation time, it is natural to examine sets of forces that 
achieve minimum propagation time.
  

\begin{defn} {\rm  Let $G=(V,E)$ be a graph and $B$ a zero forcing set of $G$.     For a set of forces $\clf$ of $B$, define $\clf^{(0)}=B$ and for $t \ge 0$,  $\clf^{(t+1)}$ is the set of vertices  $w$ such that the force $v\to w$ appears in $\clf$, $w\notin \bigcup_{i=0}^{t} \clf^{(i)}$, 
and $w$ is the only neighbor of $v$ not in $\bigcup_{i=0}^{t} \clf^{(i)}$.  The {\em propagation  time of $\clf$} in $G$, denoted  $\pt(G,\clf)$,   is the least $t_0$ such that $V=\bigcup_{t=0}^{t_0} \clf^{(t)}$.  
}
\end{defn}

 Let $G=(V,E)$ be a graph, let $B$ be a zero forcing set  of $G$, and let $\clf$ be a set of forces of $B$.  Clearly,  $\bigcup_{i=0}^{t} \clf^{(i)}\subseteq \bigcup_{i=0}^{t} B^{(i)}$ for all $t=0,\dots,\pt(G,B)$. 

\begin{defn} {\rm Let $G=(V,E)$ be a graph and let $B$ be a zero forcing set  of $G$.  A set of forces $\clf$ is  {\em efficient}   if $\pt(G,\clf)=\pt(G)$. 
Define \[\FT(G) = \{ \clf \, | \, \clf \text{ is an efficient set of forces of a minimum zero forcing set $B$ of } G\}.\]
}\end{defn} 


If $\clf$ is an efficient set of forces of a minimum zero forcing set $B$ of $G$, then $B$ is necessarily an efficient zero forcing set.   
However, not every efficient set of forces conforms to the propagation process.

   \begin{figure}[!h] 
\begin{center}\scalebox{.4}{\includegraphics{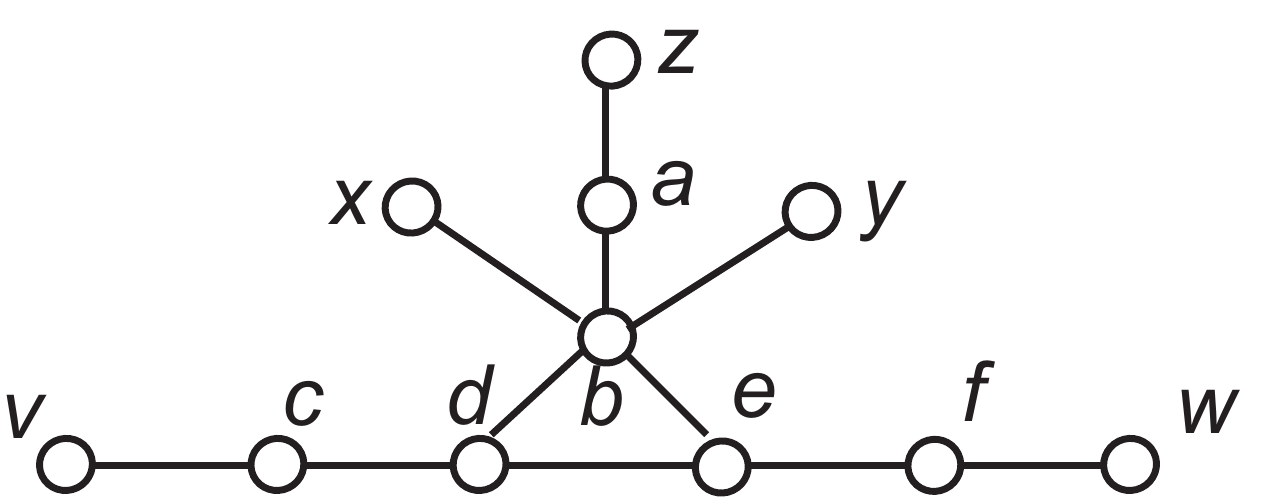}}
\caption{The graph $G$ for Example \ref{exfastnotprop}} \label{fastnotprop}
\end{center}\vspace{-5mm}
\end{figure} 

\begin{ex}\label{exfastnotprop} {\rm Let $G$ be the graph in Figure \ref{fastnotprop}.  Since every degree one vertex must be an endpoint of a maximal forcing chain and since $B=\{x,z,v\}$ is a zero forcing set, $\ZFN(G)=3$.  Since $(v,c,d,e,f,w)$ or $(w,f,e,d,c,v)$ must be a maximal  forcing chain for any set of forces of a minimum zero forcing set, $\pt(G)=5$.  Then $B$ is an efficient zero forcing set with efficient set of forces $\clf=\{v\to c, z\to a, c\to d,  a\to b, d\to e,  b\to y, e\to f, f\to w\}.$ Observe that  $b\in\clf^{(2)}$ (i.e.,  $b$ does not turn black until step $t=2$ in $\clf$), but  $b\in B^{(1)}$ ($b$ can be forced by $x$ at step $t=1$).}\end{ex}

\begin{defn} {\rm Let $G=(V,E)$ be a graph, $B$ a zero forcing set of $G$, and $\clf$ a  set of forces of $B$. 
Define $Q_0(\clf)=\T(\clf)$ and for $t = 1,\ldots, \pt(G,\clf)$, define \[Q_t(\clf) = \{ v \in V \, | \, \exists w\ \in \clf^{(\pt(G,\clf)-t+1)} \mbox{ such that }  v\to w \}.\]   }
\end{defn}

Observe that $ V=\bigcup_{t=0}^{\pt(G,\clf)} Q_t(\clf)$.

\begin{lem}\label{reverslem} Let $G = (V,E)$ be a graph, $B$  a  zero forcing set of $G$, and  $\clf$  a  set of forces of $B$. 
Then  $Q_t(\clf)\subseteq \bigcup_{i=0}^t{\Rev(\clf)}^{(i)}$. \end{lem}
\bpf
Recall that $\Rev(\clf)$ is a  set of forces of $\T(\clf)$.  The result is established by induction on $t$.   Initially, $Q_0(\clf)=\T(\clf)=\Rev(\clf)^{(0)}$. Assume that for $0\le s \le t$, $Q_s(\clf)\subseteq \bigcup_{i=0}^s{\Rev(\clf)}^{(i)}$.  Let $v\in Q_{t+1}(\clf)$.  In $\clf$,  $v\to u$ at time $\pt(G,\clf)-t$.  In $\clf$, $u$ cannot perform a force until time $\pt(G,B)-t+1$ or later, so $u\in \bigcup_{i =0}^{t} Q_i(\clf)\subseteq\bigcup_{i=0}^t{\Rev(\clf)}^{(i)}$.  If $x\in N(u)\setminus \{v\}$ then in $\clf$ $x$ cannot perform a force before time $\pt(G,\clf)-t+1$, so $x\in \bigcup_{i =0}^{t} Q_i(\clf)\subseteq \bigcup_{i=0}^t{\Rev(\clf)}^{(i)}$.  So if $v\notin \bigcup_{i=0}^t{\Rev(\clf)}^{(i)}$, then $v\in{\Rev(\clf)}^{(t+1)}$.  Thus $v\in \bigcup_{i=0}^{t+1}{\Rev(\clf)}^{(i)}$.
\epf


\begin{cor}\label{ztreversal} Let $G = (V,E)$ be a graph,  $B$  a minimum zero forcing set of $G$, and $\clf$ a  set of forces of $B$. Then \[ \pt(G,\Rev(\clf)) \leq \pt(G,\clf).\]
\end{cor}

The next result follows  from Corollary \ref{ztreversal} and Observation \ref{obsterm}.

\begin{thm}\label{ztreversal2} Let $G = (V,E)$ be a graph,  $B$  an efficient zero forcing set of $G$, and  $\clf$  an efficient set of forces of $B$. Then  $\Rev(\clf)$ is an efficient set of forces and $\T(\clf)$ is an efficient zero forcing set.  Every efficient zero forcing set is the terminus of an efficient set of forces of an efficient zero forcing set.
\end{thm}

The next result answers Question \ref{unique} negatively.

\begin{thm}\label{nonunique} Let $G $ be a connected graph of order greater than one.  Then $| \ZT(G) | \geq 2$.  
\end{thm} 
\bpf Let $B \in \ZT(G)$ and let $\clf$ be an efficient set of forces of $B$. 
By Theorem \ref{ztreversal2}, $\T(B)\in\ZT(G)$.   Since $G$ is a connected graph of order greater than one, $B \neq \T(\clf)$. 
\epf

We now consider the intersection of efficient zero forcing sets.  
The next result is immediate from Theorem \ref{ztreversal2}.

 \begin{cor}\label{termeffint} Let $G$ be a graph.  Then $\dep \bigcap_{B \in \ZT(G)} B=\bigcap_{\clf \in \FT(G)}\T(\clf)$.
 \end{cor}

\begin{quest}\label{intersect} Is there a connected graph $G$ of order at least two and a vertex $v\in V(G)$ such that $v$ is in every efficient zero forcing set? 
\end{quest}

\begin{figure}[!ht] 
\begin{center}\scalebox{.35}{\includegraphics{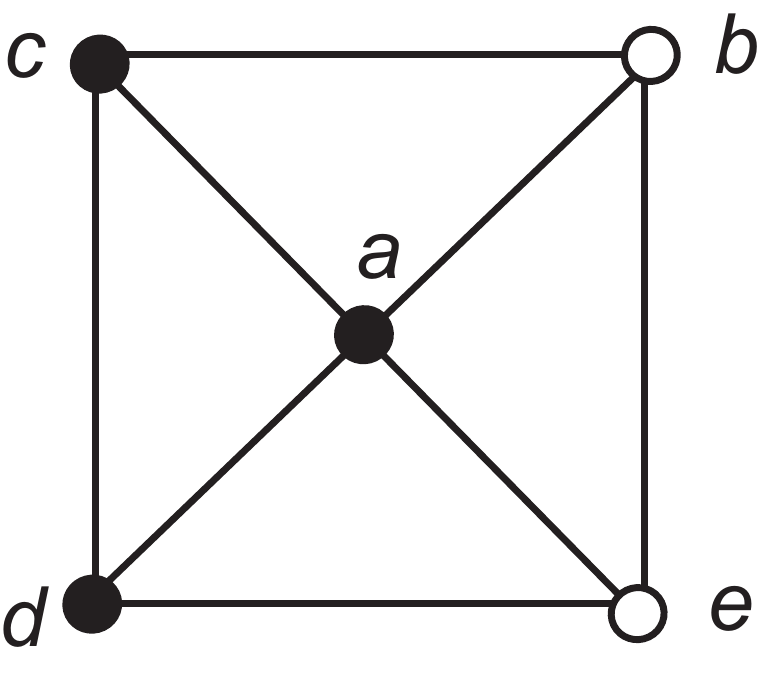}}
\caption{The wheel $W_5$}\label{wheel5}
\end{center}
\end{figure}

The next example provides an affirmative answer. 
\begin{ex}\label{W5ex}{\rm
The wheel $W_5$ is the graph shown in Figure \ref{wheel5}.  
The efficient zero forcing set $\{a,b,c\}$ of $W_5$ shows that $\pt(W_5) = 1$. Up to isomorphism, there are two types of minimum zero forcing sets in $W_5$. One set contains $a$ and  two other vertices that are adjacent to each other; the other contains three vertices other than $a$. The latter  is not an efficient zero forcing set of $W_5$, because its propagation time is $2$. The possible choices for an efficient zero forcing set are $\{a,b,c\},\{a,c,d\},\{a,d,e\},$ or $\{a,b,e\}$. Therefore, $\dep\bigcap_{B\in \ZT(G)} B=\{a\} $.
}\end{ex}

 We examine the effect of a nonforcing vertex in an efficient zero forcing set.  This result will be used in Section \ref{sslow}.
It was shown in \cite{REU09} that $\ZFN(G)=\ZFN(G-v)+1$ if and only if there exists a minimum zero forcing set $B$ containing $v$ and set of forces $\clf$ in which $v$ does not perform a force.   The proof of the next proposition is the same but with consideration restricted to an efficient zero forcing set  (the idea is that the same set of forces works for both $G$ and $G-v$, with $v$ included in the zero forcing set for $G$).

\begin{prop}\label{singletoniff} For a vertex $v$ of a graph $G$, there exists an efficient zero forcing set $B$ containing $v$ and an efficient set of forces $\clf$ in which $v$ does not perform a force if and only if $\pt(G-v) = \pt(G)$ and $\ZFN(G-v) = \ZFN(G)-1$. \end{prop}


 \section{Graphs with extreme minimum propagation time}\label{sextreme}

 For any graph $G$, it is clear that $0 \leq \pt(G)\le \PT(G) \leq \ord G-1. $  In this section we consider the extreme values $\ord G-1$,  $\ord G-2$, i.e., high propagation time, and, $0$ and $1$, i.e., low propagation time.  

\subsection{High propagation time}\label{sshigh}

The case of propagation time $\ord G-1$ is straightforward, using the well known fact \cite{Row} that $\ZFN(G)=1$ if and only if $G$ is a path.

\begin{prop} \label{ordGm1} For a graph $G$, the following are equivalent.
\ben
\item $\pt(G) =  |G| - 1$.
\item    $\PT(G) =  |G| - 1$.
\item $\ZFN(G)=1$.
\item $G$ is a path.
\een
\end{prop}

We now consider graphs $G$ that have maximum or minimum propagation time equal to $\ord G-2$.  

\begin{obs}\label{ptG2obs} For a graph $G$, 
\ben
\item $\pt(G)=\ord G-2$ implies $\PT(G) = |G| - 2$, but not conversely (see Lemma \ref{treelem} for an example).  
\item  $\pt(G) = |G| - 2$ if and only if $\ZFN(G)=2$ and   exactly one force is performed at each time for every  minimum zero forcing set. 
\item $\PT(G) = |G| - 2$  if and only if $\ZFN(G)=2$ and there exists a  minimum zero forcing set such that exactly one force performed at each time. 
\een
 \end{obs}

\begin{lem}\label{disconlem} Let $G$ be a  disconnected graph.  Then the following are equivalent.
\ben
\item $\pt(G)=\ord G -2$. 
\item $\PT(G)=\ord G -2$.
\item $G=P_{n-1} \dot{\cup} P_1$.   
\een\end{lem}
\bpf  Clearly  $G=P_{n-1} \dot{\cup} P_1\Rightarrow \pt(G)=\ord G -2\Rightarrow \PT(G)=\ord G -2$. So assume $\PT(G)=\ord G -2$.  Since $\ZFN(G)=2$, $G$ has exactly two components.  At least one component of $G$ is an isolated vertex (otherwise, more than one force occurs at time step one), and so 
$G=P_{n-1} \dot{\cup} P_1$.  \epf

A {\em path cover}  of a graph $G$  is  a set of vertex disjoint induced paths that cover all the vertices of $G$, and the path cover number $\PC(G)$ is the minimum number of paths in a path cover of $G$. It is known \cite{cancun} that for a given zero forcing set  and set of forces, the set of  maximal zero forcing chains forms a path cover and thus $\PC(G)\le \ZFN(G)$ 

A graph $G$ is a {\em graph on two parallel paths} if  $V(G)$ can be partitioned into disjoint subsets $U_1$ and $U_2$ so that the induced subgraphs $P_i=G[U_i], i=1,2$ are paths, $G$ can be drawn in the plane with the paths $P_1$ and $P_2$ as parallel line segments, and  edges between the two paths (drawn as line segments, not curves) do not cross; such a drawing is called a {\em standard drawing}.  The paths $P_1$ and $P_2$ are called the {\em parallel paths} (for this representation of $G$ as a graph on two parallel paths).  

  Let $G$ be a graph on two parallel paths  $P_1$ and $P_2$.
  If $v\in V(G)$, then $\pth(v)$ denotes the parallel path that contains $v$ and $\npath(v)$ denotes the other of the parallel paths. Fix an ordering of the vertices in each of $P_1$ and $P_2$ that is increasing in the same direction for both paths in a standard drawing.  With this ordering, let $\first(P_i)$ and $\last(P_i)$ denote the first and last vertices of $P_i, i=1,2$.   If $v,w\in V(P_i)$, then $v\prec w$ means $v$ precedes $w$ in the order on $P_i$.  Furthermore, if $v\in V(P_i)$ and $v\ne \last(P_i)$, $\next(v)$ is the neighbor of $v$ in  $P_i$ such that  $v\prec \next(v)$; $\prev(v)$ is defined analogously (for $v\ne \first(P_i)$). 

Row \cite{Row} has shown that $Z(G)=2$ if and only if $G$ is a graph on two parallel paths. Observe that for any graph having $\ZFN(G)=2$,  a  set of forces $\clf$ of a minimum zero forcing set  naturally produces a representation of $G$ as a graph on two parallel paths with the parallel paths being the maximal forcing chains. The ordering of the vertices in the parallel paths is  the forcing order.

\begin{lem}\label{treelem} For a tree $G$, $\PT(G)=\ord G -2$ if and only if $G=S(1,1,n-3)$ (sometimes called a T-shaped tree). The graph $K_{1,3}$ is  the only tree for which $\pt(G)=\ord G -2$. 
\end{lem}
\bpf  
It is clear that  $\PT(S(1,1,n-3))=n-2$, and $\pt(K_{1,3})=2$. 

Suppose first that $G$ is a tree such that $\PT(G)=\ord G -2$.  Then $G$ is a graph on two parallel paths $P_1$ and $P_2$.     There is exactly one edge $e$ between the two paths.  Observe that $e$ must have an endpoint not in $\{\first(P_i), \last(P_i), i=1,2\}$, so without loss of generality $\first(P_1)\ne\last(P_1)$ and neither $\first(P_1)$ nor $\last(P_1)$ is an endpoint of $e$. 
If $G$ is a graph with multiple vertices in each of $P_1, P_2$ (i.e., if $\first(P_2)\ne \last(P_2)$), then no matter which minimum zero forcing set we choose, more than one force will occur at some time. So assume $V(P_2)$ consists of a single vertex $w$.  If the parallel paths were constructed from a minimum zero forcing set $B$, then $w\in B$, and without loss of generality $B=\{\first(P_1),w\}$.  If $N(w)\ne N(\first(P_1))$, then  at time one,  two vertices would be forced.  Thus $N(w)=N(\first(P_1))$ and  $G=S(1,1,n-3)$.  

 Now suppose that $G$ is a tree such that $\pt(G)=\ord G -2$.   This implies $\PT(G)=\ord G -2$, so  $G=S(1,1,n-3)$.  Since $n-3>1$ implies $\pt(S(1,1,n-3))<n-2$, $G=S(1,1,1)=K_{1,3}$.
  \epf



\begin{figure}[!ht] 
\begin{center}\scalebox{.45}{\includegraphics{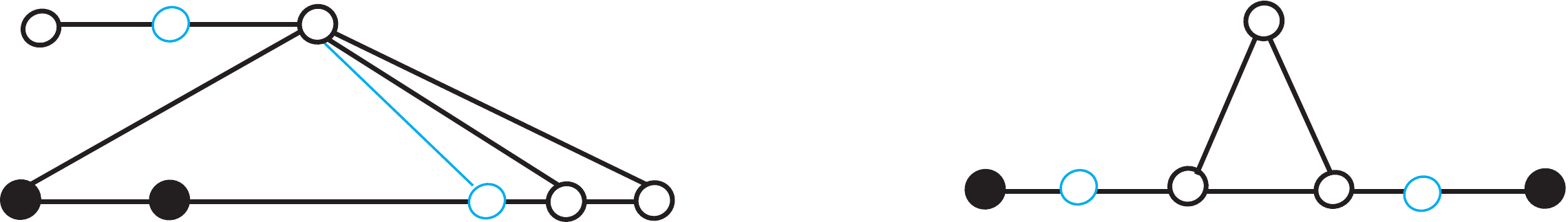}}
\caption{Graphs $G$ and minimum zero forcing sets $B$ such that  $\pt(G,B)<\ord G -2 $ (where {\blue light} vertices may be absent or repeated and similarly for {\blue light} edges)} \label{fzzbad}
\end{center}
\end{figure}

\begin{obs}  If $G$ is one of the graphs shown in Figure \ref{fzzbad}, then $\pt(G)<\ord G-2$, because the black vertices are a minimum zero forcing set $B$ with $\pt(G,B)<\ord G -2$.
\end{obs}

For any graph and vertices $x,y$,  $x\sim y$ denotes that $x$ and $y$ are adjacent, and  $xy$ denotes the edge with endpoints $x$ and $y$.

 \begin{defn}{\rm  
A graph $G$ on two parallel paths $P_1$ and $P_2$ is a {\em zigzag graph} if it satisfies the following conditions:
\ben
\item  There is a path $Q=(z_1,z_2,\dots,z_\ell)$ that alternates between the two paths $P_1$ and $P_2$ such that:
\ben
\item $z_{2i-1}\in V(P_1)$ and $z_{2i}\in V(P_2)$ for $i=1,\dots,\lfloor\frac{\ell+1}2\rfloor$;
\item $z_j\prec z_{j+2}$ for $j=1,\dots,\ell-2$.
\een
\item Every edge of $G$ is an edge of one of  $P_1, P_2$, or $Q$, or is an edge of the form
\[z_jw \mbox{ where } 1<j<\ell, w\in\npath(z_j), \mbox{ and } z_{j-1}\prec w\prec z_{j+1}.\]
\een
The number $\ell$ of vertices in $Q$ is called the {\em zigzag order}.}
\end{defn}

\begin{figure}[!ht] 
\begin{center}\scalebox{.45}{\includegraphics{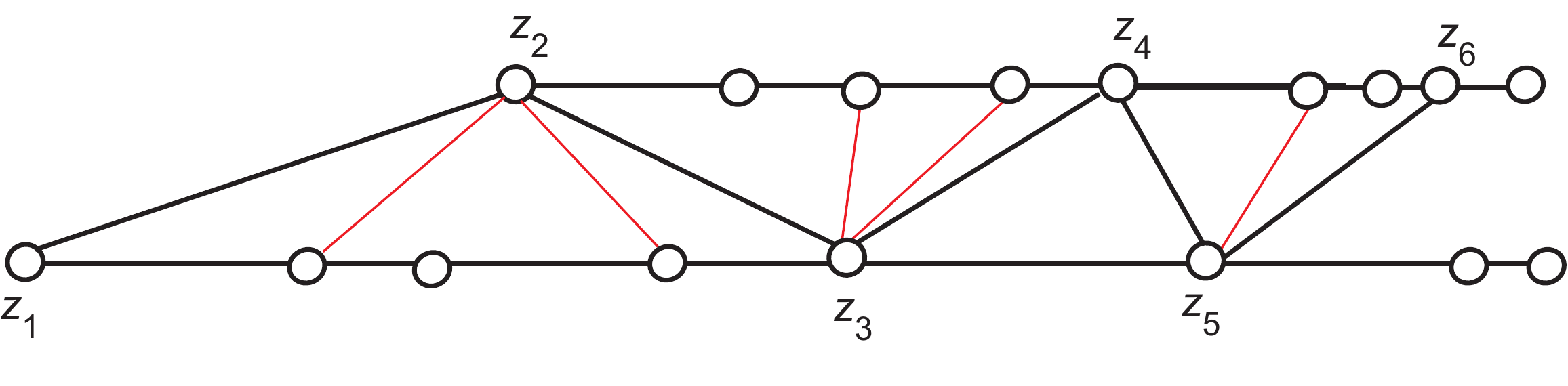}}
\caption{A  zigzag graph (with $P_1, P_2$ and $Q$ in black)} \label{zigzagfig1}
\end{center}
\end{figure}

\begin{figure}[!ht] 
\begin{center}\scalebox{.45}{\includegraphics{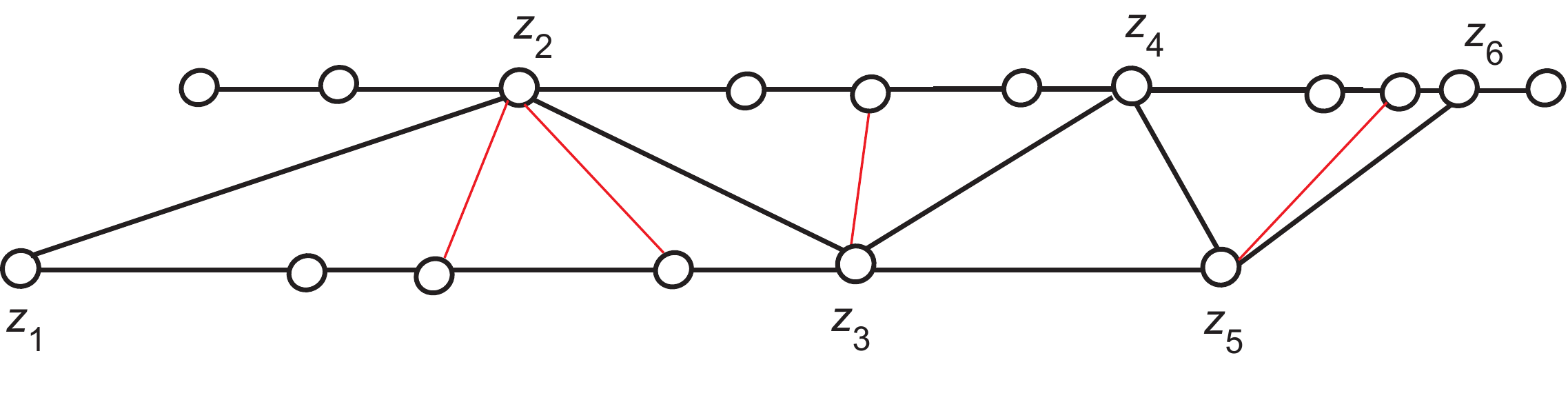}}
\caption{A zigzag graph $G$  (with $P_1, P_2$ and $Q$ in black) having $\pt(G)=\ord G-2$ } \label{fzzpt1}
\end{center}
\end{figure}

 Examples of zigzag graphs  are shown in Figures \ref{fzzbad},  \ref{zigzagfig1}, and \ref{fzzpt1}.   

\begin{thm} \label{ordGm2pt} Let $G$ be a graph. Then $\pt(G) =  |G| - 2$ if and only if $G$ is one of the following:
\ben
\item $P_{n-1}\dot{\cup} P_1$.
\item $K_{1,3}$.

\item\label{speczz} A 
zigzag graph of zigzag order $\ell$ such that all of the following conditions are satisfied:
\ben
 \item\label{ecx1} $G$ is not isomorphic to one of the graphs shown in Figure \ref{fzzbad}.
 \item\label{ecx2} $\deg (\first(P_1))>1$ or $\deg (\first(P_2))>1$ (both paths cannot begin with degree-one vertices).
  \item\label{ecx3} $\deg (\last(P_1))>1$ or $\deg (\last(P_2))>1$   (both paths cannot end with degree-one vertices).
   \item\label{ecx4} $z_2\ne \first(P_2)$ or $z_2\sim \next(z_1)$ 
 \item\label{ecx5} $z_{\ell-1}\ne \last(\pth(z_{\ell-1}))$ or $z_{\ell-1}\sim \prev(z_\ell)$ 
 \een \een  
\end{thm}

An example of a zigzag graph satisfying  conditions (\ref{ecx1}) --  (\ref{ecx5}) is shown in Figure \ref{fzzpt1}.

\bpf  Assume   $\pt(G) = |G| - 2$.   If $G$ is disconnected or a tree, then  $G$ is $P_{n-1}\dot{\cup} P_1$ or $K_{1,3}$ by Lemmas \ref{disconlem} and \ref{treelem}.  So assume $G$ is connected and has a cycle.  

First we identify $P_1, P_2$ and $Q$ for $G$:  By Observation \ref{ptG2obs}, there exists a minimum zero forcing set $B$ of cardinality 2 such that  exactly one force is performed at each time for $B$.  
Renumber the vertices of $G$ as follows:  vertices $V(G)=\{-1,0,1,2,\dots,n-2\}$, zero forcing set  $B=\{-1,0\}$ with $0\to 1$, and vertex $t$ is forced at time $t$.  
Then $G$ is a graph on two parallel paths  $P_1$ and $P_2$, which are the two maximal forcing chains (with the path order being the forcing order).  Observe that $\deg (0)\le 2$ and $\deg (-1)\ge 2$, because $0$ can immediately force and $-1$ cannot.  
If $\deg (-1)=2$ and $\ord G>3$, then choose $P_1$ to be $\pth(-1)$, and let $z_1=-1, \ z_2=\max N(-1) \cap P_2$.   Otherwise, choose  $P_1$ to be $\pth(0)$ and let $z_1=\min N(-1), \ z_2=-1$.     
For $j\ge 2$, define $z_{j+1}=\max N(z_j) \cap \npath(z_j)$ until $N(z_j) \cap \npath(z_j)=\emptyset$.   Define $Q=(z_1,\dots,z_\ell)$.  With this labeling, $G$ is a zigzag graph. 

Now we show that $G$  satisfies conditions (\ref{ecx1}) --  (\ref{ecx5}).  Since $\pt(G)=\ord G-2$, $G$ is not isomorphic to one of the graphs shown in Figure \ref{fzzbad}, i.e., condition (\ref{ecx1}) is satisfied.  Since $-1$ is the first vertex  in one of the paths and $\deg(-1)\ge 2$, condition (\ref{ecx2}) is satisfied.  The remaining conditions must be satisfied or there is a different zero forcing set of two vertices with lower propagation time: if (\ref{ecx3}) fails, use $B=\{\last(P_1), \last(P_2)\}$; if (\ref{ecx4}) fails, then $z_2= \first(P_2)$ and $z_2\not\sim \next(z_1)$, so use $B=\{\first(P_1), \next(z_1)\}$; if (\ref{ecx5}) fails, this is analogous to (\ref{ecx4}) failing, so use $B=\{\last(\pth(z_\ell)), \prev(z_{\ell})\}$.

For the converse,   $\pt(G)=\ord G -2$ for  $G=P_{n-1}\dot{\cup} P_1$ or $G=K_{1,3}$ by Lemmas \ref{disconlem} and  \ref{treelem}.  So assume   $G$ is a zigzag graph satisfying conditions (\ref{ecx1}) --  (\ref{ecx5}).  The sets $B_1=\{\first(P_1),\first(P_2)\}$ and $B_2=\{\last(P_1),\last(P_2)\}$ are minimum zero forcing sets of $G$, and $\pt(G,B_i)=\ord G-2$ for $i=1,2$. 
If $z_2=\first(P_2)$, so $\first(P_2)\sim \next(z_1)$, then $B_3=\{\first(P_1),\next(z_1)\}$ is a zero forcing set and $\pt(G,B_3)=\ord G-2$.  
If $z_{\ell-1}\ne \last(\pth(z_{\ell-1}))$, so  $\last(\pth(z_{\ell-1}))\sim \prev(z_\ell)$, then $B_4=\{\last(\pth(z_\ell)),\prev(z_\ell)\}$ is a zero forcing set and $\pt(G,B_4)=\ord G-2$.  If $G$ is not isomorphic to one of the graphs shown in Figure \ref{fzzbad}, these are the only minimum zero forcing sets.  \epf

\vspace{-5mm}

\subsection{Low propagation time}\label{sslow}

\begin{obs} \label{pt0} For a graph $G$, the following are equivalent.
\ben
\item $\pt(G) =  0$.
\item    $\PT(G) = 0$.
\item $\ZFN(G)=\ord G$.
\item $G$ has no edges.
\een
\end{obs}

Next we consider $\pt(G)=1$.  From Remark \ref{ptPTbds} we see that if $\pt(G) =  1$, then $\dep \ZFN(G) \geq \frac{|G|}{2}$. 
The converse of this statement is false:

\begin{ex}\label{exMH1}{\rm  Let $G$ be the graph obtained from $K_4$ by appending a leaf to one vertex.  Then $\ZFN(G)  = 3 >|G|/2$ and $\pt(G) =  2$. }\end{ex}

\begin{thm}\label{pt1force}  Let $G=(V,E)$ be a connected graph such that  $\pt(G)=1$.  For $v\in V$,  $\dep v\in \bigcap_{B \in \ZT(G)} B$ if and only if for every $B\in\ZT(G)$, $\ord{B^{(1)}\cap N(v)}\ge 2$.
\end{thm}
\bpf  
If for every $B\in\ZT(G)$,  $\ord{B^{(1)}\cap N(v)}\ge 2$, then $v$ cannot perform a force in an efficient set of forces, so $v\in\T(\clf)$ for every $\clf\in\FT(G)$.  Thus $v\in \bigcap_{\clf \in \FT(G)}\T(\clf)=\bigcap_{B \in \ZT(G)} B$.

Now suppose $v\in \bigcap_{B \in \ZT(G)} B$ and let $ B\in \ZT(G)$.  If $v$ performs a force in an efficient set of forces  $\clf$ of $B$, then  $v\notin\T(\clf)$.  By Theorem \ref{ztreversal2}, $\T(\clf)\in \ZT(G)$, so $v$ cannot perform a force in any such $\clf$.  Since $v$ cannot perform a force, $\ord{B^{(1)}\cap N(v)}\ne 1$.  It is shown in \cite{Zplus} that (assuming the graph is connected and of order greater than one) every vertex of a minimum zero forcing set   must have a neighbor not in the zero forcing set.  Since $\pt(G)=1$, $\ord G>1$, and so $\ord{B^{(1)}\cap N(v)}\ge 2$.  \epf

We now consider the case of a graph $G$ that has $\pt(G)=1$ and $\ZFN(G)=\frac 1 2 \ord G$.  Examples of such graphs include the hypercubes $Q_s$ \cite{AIM}.     

\begin{defn}{\rm Suppose $H_1=(V_1,E_1)$ and $H_2=(V_2,E_2)$ are graphs of equal order and $\mu:V_1\to V_2$ is a bijection. Define the {\em matching graph} $(H_1,H_2,\mu)$ to be the graph constructed as the disjoint union of $H_1, H_2$ and the perfect matching between $V_1$ and $V_2$ defined by $\mu$.}
\end{defn} 

Matching graphs play a central role in the study of graphs that have propagation time one.  

\begin{prop}\label{matchprop} Let $G=(V,E)$ be a graph.  Then any two of the following conditions imply the third.
\ben
\item\label{cZ} $\ord G = 2\ZFN(G)$.  
\item\label{cpt} $\pt(G) = 1$.
\item\label{cmatch} $G$ is a matching graph
\een
\end{prop}
\bpf $(\ref{cZ}) \ \& \ (\ref{cpt})  \Rightarrow (\ref{cmatch})$: Let $B$ be an efficient zero forcing set of $G$ and let $\overline{B}=V\setminus B$.   Since $\ord B=\frac 1 2 \ord G$ and $\pt(G)=1$, every element  $b \in B$ must perform a force at time one.    Thus  $\ord{N(b) \cap \overline{B}}=1$ and there exists a perfect matching between $B$ and $\overline{B}$ defined by $\mu:B \to {\overline B}$ where  $\mu(b)\in N(b) \cap \overline{B}$.   Then $G=(B,{\overline B},\mu)$.

For the remaining two parts, assume $G=(H_1,H_2,\mu)$ and $n=\frac 1 2 \ord G$ ($=\ord{H_1}=\ord{H_2}$).

$(\ref{cZ}) \ \& \ (\ref{cmatch})  \Rightarrow (\ref{cpt})$:  Since $\ZFN(G)=n$, $H_1$ is a minimum zero forcing set and $\pt(G,H_1)=1$.

$(\ref{cpt}) \ \& \ (\ref{cmatch})  \Rightarrow (\ref{cZ})$:   Since $\pt(G)=1$, $\ZFN(G)\ge n$, and $\ZFN(G)\le n$  because $H_1$ is a zero forcing set with $\pt(G,H_1)=1$.
\epf

We examine conditions that ensure $\ZFN((H_1,H_2,\mu))=\ord{H_i}$ and thus $\pt((H_1,H_2,\mu))=1$.
The choice of matching $\mu$ affects the zero forcing number and propagation time, as the next two examples show.

\begin{ex}\label{C5P2}{\rm  The  Cartesian product $C_5\Box P_2$ is $(C_5,C_5,\iota)$, where $\iota$ is the identity mapping. It is known \cite{AIM} that $\ZFN(C_5\Box P_2)  = 4$   and thus $\pt(C_5\Box P_2)  >1$. }\end{ex}

\begin{ex}\label{Pet}{\rm  The Petersen graph $P$ can be constructed as $(C_5,C_5,\mu_P)$ where $\mu_P=\left(\begin{array}{rrrrr}
1 & 2 & 3 & 4 & 5 \\
1 & 4 & 2 & 5 & 3
\end{array}\right)$. It is known \cite{AIM} that $\ZFN(P)  = 5$   and thus $\pt(P)  =1$. }\end{ex}

Let $c(G)$ denote the number of components of $G$.

\begin{thm}\label{H1H2mucomponent} Let $\ord {H_1}=\ord {H_2}=n$ and let $\mu:H_1\to H_2$ be a bijection.  If $\pt((H_1,H_2,\mu))=1$, then $c(H_1)=c(H_2)=c((H_1,H_2,\mu))$.
\end{thm}
\bpf  
Assume it is not the case that $c(H_1)=c(H_2)=c((H_1,H_2,\mu))$.  This implies  $\mu$ is not the union of perfect matchings between the components of $H_1$ and the components of $H_2$. Without loss of generality, there is a component $H_1[C_1]$ of $H_1$ that is not matched within a single component of $H_2$. Then there exist vertices $u$ and $v$ in $C_1$ such that $\mu(v) \in C_v$,  $\mu(u) \in C_u$, and $H_2[C_v]$ and $H_2[C_u]$ are separate components of $H_2$. We show that there is a zero forcing set of size $n-1$ for $(H_1,H_2,\mu)$, and thus pt$((H_1,H_2,\mu))>1$. Let $B_1 = C_1 \backslash (\mu^{-1}(C_v) \bigcup \{u\})$, $B_2 = V_2 \backslash (\mu(B_1) \bigcup \mu(u))$, and $B = B_1 \bigcup B_2$, so $|B|=n-1$. Then 1) $x \rightarrow \mu^{-1}(x)$ for $x \in C_v$, 2) $v \rightarrow u$, 3) $y \rightarrow \mu(y)$ for $y \in C_1 \backslash \mu^{-1}(C_v)$, and 4) $z \rightarrow \mu^{-1}(z)$ for all $\mu^{-1}(z)$ in the remaining components of $H_1$. Therefore B is a zero forcing set,  $\ZFN((H_1,H_2,\mu))\le n-1$, and thus $\pt((H_1,H_2,\mu))>1$.
\epf

\begin{thm}\label{KrmatchH}  Let $\ord H=n$ and let $\mu$ be a bijection of vertices of   $H$ and $K_n$  (with $\mu$ acting on the vertices of $H$).  Then   $\pt((H,K_n,\mu))=1$ if and only if $H$ is connected.  
\end{thm}
\bpf  
If $H$ is not connected, then $\pt((H,K_n,\mu))\ne 1$ by Theorem \ref{H1H2mucomponent}. Now assume $H$ is connected and let $G=(H,K_n,\mu)$. Let $B\subseteq V(G)$ with $\ord B =n-1$.  We show $B$ is not a zero forcing set. This implies $\ZFN(G)=n$ and thus $\pt(G)=1$.  Let $X=V(K_n)$ and $Y=V(H)$.    For $x\in X$, $x$ cannot perform a force until at least $n-1$ vertices in $X$ are black. If $\ord {X\cap B} =n-1$ then no force can be performed.  So assume $\ord {X\cap B} \le n-2$. Until at least $n-1$ vertices in $X$ are black, all forces must be performed by vertices in $Y$.   We show that no more than $n-2$ vertices in $X$ can turn black.  Perform all forces of the type $y\to y'$ with $y,y'\in Y$. For each such force, $\mu(y)$ must be black already. Thus at most $\ord{X\cap B}$ such forces within $Y$ can be performed.   So there are now at most $\ord {Y\cap B}+\ord{X\cap B}=n-1$ black vertices in $Y$.  Note first that if at most $n-2$ vertices of $Y$ are black, then after all possible forces from $Y$ to $X$ are done, no further forces are possible, and at most $n-2$ vertices in $X$ are black.  So assume $n-1$ vertices of $Y$ are black.  Let $w\in Y$ be white.  Since $H$ is connected, there must be a neighbor $u$ of $w$ in $Y$, and $u$ is black.  Since $u\in N(w)$ and $w$ is white, $u$ has not performed a force.  If $\mu(u)$ were black,  there would be at most $n-2$ black vertices in $Y$, so $\mu(u)$ is white.    After preforming all possible forces from $Y$ to $X$, at most $n-2$ vertices in $X$ are black because all originally black vertices $x$ have $\mu^{-1}(x)$ black, there are $n-1$ black vertices in $Y$, and  $u$ cannot perform a force at this time (since both $w$ and $\mu(u)$ are white).  Thus not more than $n-2$ vertices of $X$ can be forced, and $B$ is not a zero forcing set.
 \epf

The Cartesian product of $G$ with $P_2$ is one way of constructing matching graphs, because $G \, \square \, P_2=(G,G,\iota)$.  
 Examples of graphs $G$ having $\ZFN(G \, \square \, P_2) = |G|$ include the complete graph $K_r$ and hypercube $Q_s$ \cite{AIM}.  Since $\ZFN(G\,\Box\, P_2)\le 2\ZFN(G)$ \cite{AIM}, to have $\ZFN(G\,\Box\, P_2)= \ord G$ it is necessary that $\ZFN(G)\ge \frac{\ord G} 2$, but that condition is not sufficient.

\begin{ex} \label{ex2} {\rm Observe that $\ZFN(K_{1,r})=r-1\ge \frac 1 2 \ord{K_{1,r}}$ for $r\ge 3$. 
But  $\ZFN(K_{1,r} \, \square \, P_2) = r < \ord{K_{1,r}}$, so $\pt(K_{1,r}\, \square \, P_2) \geq 2$. }\end{ex}

The next theorem provides conditions that ensure that  iterating the Cartesian product with $P_2$ gives a graph with propagation time one.   Recall that one of the original motivations for defining the zero forcing number was to bound maximum nullity, and the interplay between these two parameters is central to the proof of the next theorem.  Let $G=(\{v_1,\dots,v_n\},E)$ be a graph. The {\em set of symmetric matrices described by  $G$} is
$\SG=\{A\in\Rnn : A^T=A \mbox{ and for } i\ne j, a_{ij} \ne 0 \Leftrightarrow v_iv_j\in E\}.$
The {\em maximum nullity of $G$} is 
$\M(G)=\max\{{\rm null}\, A: A\in \SG\}.$  It is well known  \cite{AIM} that $\M(G)\le\ZFN(G)$.
The next theorem provides conditions that are sufficient to iterate the construction of taking the Cartesian product of a graph and $P_2$ and obtain minimum propagation time equal to one. 

\begin{thm}\label{GboxP2}  Suppose $G$ is a graph with $\ord G=n$ and there exists a matrix $L\in\sym(G)$ such that $L^2=I_n$.  Then \[\M(G\,\Box\, P_2)= \ZFN(G\,\Box\, P_2)=n \mbox{ and } \pt(G\,\Box\, P_2)=1.\]  Furthermore, for \[\hat L=\frac 1 {\sqrt 2}\mtx{L & I_n\\I_n & -L}\] $\hat L\in \sym(G\,\Box\, P_2)$ and ${\hat L}^2=I_{2n}$.
\end{thm}
\bpf
Given the $n\x n$ matrix $L$, define
\[
H=\mtx{L&I_n\\ I_n&L}.
\]
Then $H,\hat L\in\sym(G\,\Box\, P_2)$ and ${\hat L}^2=I_{2n}$. Since $\mtx{I_n&0\\ -L & I_n} \mtx{L&I_n\\ I_n&L}=\mtx{L&I_n\\ 0&0}$, $\nul(H)=n$.
Therefore,  $\M(G\,\Box\, P_2)\ge n$.  Then 
\[n\le \M(G\,\Box\, P_2)\le\ZFN(G\,\Box\, P_2)\le n\]
so we have equality throughout.  Since $G\,\Box \,P_2$ is a matching graph, $\pt(G\,\Box \,P_2)=1$ by Proposition \ref{matchprop}.
\epf

Let $G\,(\Box\, P_2)^s$ denote the graph constructed by starting with $G$ and performing the Cartesian product with $P_2$ $s$ times.  For example, the hypercube $Q_s= P_2\,(\Box\, P_2)^{s-1}$, and the proof given in \cite{AIM} that $\M(Q_s)=\ZFN(Q_s)=2^{s-1}$ is the same as the proof of Theorem \ref{GboxP2} using  the matrix $L=\mtx{0 & 1\\1 & 0}\in\sym(P_2)$.

\begin{cor}\label{GboxP2s}  Suppose $G$ is a graph such that there exists a matrix $L\in\sym(G)$ such that $L^2=I_{\ord G}$.  Then for $s\ge 1$,  \[\M(G\,(\Box\, P_2)^s)=\ZFN(G\,(\Box\, P_2)^s)=\ord G 2^{s-1} \mbox{ and }\pt(G\,(\Box\, P_2)^s)=1.\]  
\end{cor}

\begin{cor}\label{KrboxP2s}  For $s\ge 1$, $\M(K_n\,(\Box\, P_2)^s)=\ZFN(K_n\,(\Box\, P_2)^s)=n 2^{s-1}$ and $\pt(K_n\,(\Box\, P_2)^s)=1$.  
\end{cor}
\bpf  Let $L=I_n-\frac 2 n J_n$, where $J_n$ is the $n\x n$ matrix having all entries equal to one.  Then $L\in\sym(K_n)$ and $L^2=I_n$. 
\epf

Observe that if the matrix $L$ in the hypothesis of Theorem \ref{GboxP2} is symmetric (and thus is an orthogonal matrix), then the matrix $\hat L$ in the conclusion also has these properties.  The same is true for  Corollary \ref{KrboxP2s}, where $L=I_n-\frac 2 n J_n$ is a Householder transformation.

 We have established a number of constructions that provide matching graphs having propagation time one.  For example,   $\pt(P_2\,(\Box\, P_2)^s)=1$, $\pt(K_n\,(\Box\, P_2)^s)=1$, and if $H$ is connected, then $\pt((K_n,H,\mu))=1$ for every matching $\mu$.  But the general question remains open.
 
 \begin{quest} Characterize matching graphs $(H_1,H_2,\mu)$ such that $\pt((H_1,H_2,\mu))=1$.
 \end{quest}

We can  investigate when $\pt(G)=1$ by deleting vertices that are in an efficient zero forcing set but do not perform a force in an efficient set of forces.  
The next result is a consequence of Proposition \ref{singletoniff}.

\begin{cor}\label{cor2}   Let $G$ be a graph with $\pt(G) = 1$, $B$  an efficient zero forcing set of $G$ containing $v$, and $\clf$ an efficient set of forces of $B$ in which $v$ does not perform a force. Then $\pt(G-v) = \pt(G) = 1$. \end{cor}

\begin{defn}{\rm Let $G$ be a graph with $\pt(G) = 1$, $B$  an efficient zero forcing set of $G$, $\clf$ an efficient set of forces of $B$, and $S$  the set of vertices in $B$ that do not perform a force. Define $V'=V\setminus S$, $G'=G[V']=G- S$,  $B'=B\setminus S$, and $\overline{B'}=V'\setminus B'$. The graph $G'$ is called a {\em prime} subgraph of $G$ with associated zero forcing set $B'$.
}\end{defn}
\begin{obs} \label{pt1structure} Let $G$ be a graph with $\pt(G) = 1$.   For the prime subgraph $G'$ and associated zero forcing set $B'$ defined from    an efficient zero forcing set $B$ and  efficient set of forces $\clf$ of $B$:
 \ben
 \item $\overline{B'}=V\setminus B$.
\item $\ord B' = \ord  {\overline{B'}}$ and   $|G'| = 2|B'|$.
\item   $G'$ is the matching graph defined by  $G[B'], G[\overline{B'}]$ 
and  $\mu:B' \to \overline{B'}$ defined by $\mu(b)\in(N(b) \cap \overline{B'})$.  
\item $B'$ and  ${\overline B'}$ are efficient zero forcing sets of $ G'$.
\item $\pt(G') = 1$.
\een
\end{obs}

 It is clear that if $G=(V,E)$ has no isolated vertices, $\pt(G)=1$, and if $\hat G$ is constructed from $G$ by  adjoining a new vertex $v$ adjacent to every $u\in V(G)$, then $\pt(\hat G)=1$.

We now return to considering  \ $\bigcap_{B \in \ZT(G)} B$, specifically in the case of propagation time one. 
We have a corollary of Theorem \ref{pt1force}. 

\begin{cor}\label{deg4ex2} Let $G = (V,E)$ be a graph such that $\pt(G) = 1$. If $\dep v \in \bigcap_{B \in \ZT(G)} B$, then $\deg v \geq 4$. \end{cor}
 \bpf Let $\dep v \in \bigcap_{B \in \ZT(G)} B$.  Since $\pt(G)=1$, for any efficient set $\clf$ of $B$, $\T(B)=B^{(1)}$.  By Theorem \ref{pt1force}, $\ord{B^{(1)}\cap N(v)}\ge 2$, so $\ord{\T(B)\cap N(v)}\ge 2$. Since $B=\T(\Rev(\clf))$,   $\ord{{B}\cap N(v)}\ge 2$.  Thus $\deg v\ge 4$. \epf

Note that Corollary \ref{deg4ex2} is false without the hypothesis that $\pt(G)=1$, as the next example shows.

\begin{ex} \label{deg4ex3}{\rm  For   the graph $G$ in Figure \ref{deg42}, $\ZFN(G) = 3$. Every minimum zero forcing set   $B$ must contain one of $\{a,c\}$ and one of $\{x,z\}$; without loss of generality, $a,x\in B$. If $v\in B$ then $\pt(G,B)=2$; if not then $c$ or $z$ is in $B$ and $\pt(G,B)=3$. Thus $v$ is in every efficient zero forcing set.

  \begin{figure}[!ht] 
\begin{center}\scalebox{.4}{\includegraphics{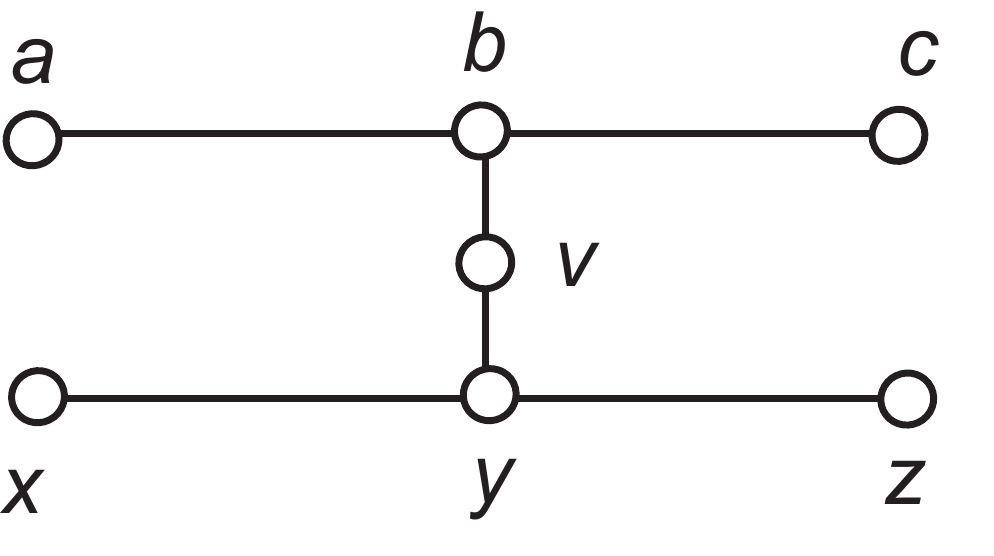}}
\caption{A graph $G$ with  $v\in\bigcap_{B \in \ZT(G)} B$ and $\deg v <4$ } \label{deg42}
\end{center}
\end{figure}
}\end{ex} \vspace{-5mm}

\begin{prop}\label{propintersect} Let  $G=(V,E)$ be a graph and $v\in V$.  If $\deg v> \ZFN(G)$ and $\pt(G)=1$, then  $v\in \bigcap_{B \in \ZT(G)} B$.
\end{prop}  
\bpf   Suppose $v\notin B\in\ZT(G)$, and let $\clf$ be an efficient set of forces of $B$.   Then $v$ performs a force in the efficient set $\Rev(\clf)$ of $\T(\clf)$.  Since every force is performed at time 1, $\deg v\le \ZFN(G)$.  \epf

The converse of Proposition \ref{propintersect} is false, as the next example demonstrates.

\begin{ex} \label{deg4ex}{\rm  Let $G$ be the graph  in Figure \ref{deg4}. 
  \begin{figure}[!ht] 
\begin{center}\scalebox{.8}{\includegraphics{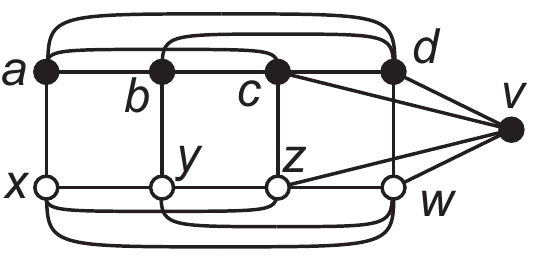}}
\caption{A graph $G$ with  $v\in\bigcap_{B \in \ZT(G)} B$ and $\deg v <\ZFN(G)$ } \label{deg4}
\end{center}
\end{figure}
It can be verified that $\ZFN(G) = 5$. Then $B_1=\{a,b,c,d,v\}$ and $B_2=\{x,y,z,w,v\}$ are  minimum zero forcing sets and $\pt(G,B_1)=\pt(G,B_2)=1$, so $\pt(G) = 1$. Let $B$ be a zero forcing set  of $G$ not containing vertex $v$.  In order to have $\pt(G,B)=1$, some neighbor of $v$ must be able to force $v$ immediately. Without loss of generality, this neighbor is $d$.  Then $a,b,c,d,w\in B$.  The set $\{a,b,c,d,w\}$ is a minimum zero forcing set but has propagation time  $2$, because no vertex can force $z$ immediately.  Thus $\bigcap_{B \in \ZT(G)} B = \{v\}$, and observe that $\deg v = 4 < 5 =\ZFN(G)$.
}\end{ex}


\section{Relationship of propagation time and diameter}\label{sdiam}

In general, the diameter and the propagation time of a graph are not comparable. Let $G$ be the dart (shown in Figure \ref{dart}). Then $\diam(G) = 2 < \pt(G) = 3$. On the other hand, $\diam(C_4) = 2 > 1 = \pt(C_4)$.

Although it is not possible to a  obtain a direct ordering relationship between  diameter and propagation time in an arbitrary graph, diameter serves as an upper bound for propagation time in the family of trees. To demonstrate this, we need some definitions. 
 The {\em walk} $v_1v_2\dots v_p$ in $G$ is the subgraph with vertex set $\{v_1,v_2,\dots, v_p\}$  and edge set $\{v_1v_2, v_2v_3, \dots, v_{p-1}v_p\}$ (vertices and/or edges may be repeated in these lists but are not repeated in the vertex and edge sets). A  {\em trail} is a walk with no repeated edges (vertices may be repeated; 
a {path} is a trail with no repeated vertices). The {\em length} of a trail $P$, denoted by $\len (P)$, is the number of edges in $P$.  We show in Lemma \ref{ztraillem} below that for any graph $G$ and minimum zero forcing set $B$, there is a trail of length at least $\pt(G,B)$.  A trail produced by the method in the proof is illustrated in the Example \ref{extrail} below.


\begin{lem}\label{ztraillem} Let $G$ be a graph and let $B$ be a minimum zero forcing set of $G$. Then there exists a trail $P$ such that $\pt(G,B)\le\len(P)$.
\end{lem}
\bpf
Observe that if $u,v \in V(G)$ such that $u$ forces $v$ at time $t>1$, then $u$ cannot force $v$ at time $t-1$.  Thus  either $u$ was forced at time $t-1$ or some neighbor of $u$ was forced at time $t-1$. So there is a path $w u v$, where $w$ forces $u$ at time $t-1$, or a path $w x u v$, where $w$ forces $x$ at time $t-1$ and $x$ is a neighbor of $u$.

 We construct a trail $v_{-p} v_{-p+1} v_{-p+2} \dots v_0$, such that for each time $t$, $1 \le t \le \pt(G,B)$, there exists an $i_t$, $-p \le i_t \le -1$, such that $v_{i_t}$ forces $v_{i_t+1}$ at time $t$.
 Begin with $t=\pt(G,B)$ and work backwards to $t=1$ to construct the trail,  using negative numbering.  To start,  there is some vertex $v_0$ that is forced by a vertex   $v_{-1}$ at time $t=\pt(G,B)$; the  trail is now $v_{-1}v_0$. Assume the trail $v_{-j}\dots v_0$ has been constructed so that for each time $t=\ell, \dots, \pt(G,B)$, there exists an $i_t$, such that $v_{i_t}$ forces $v_{i_t+1}$ at time $t$. If $\ell>1$, then $v_{-j}\to v_{-j+1}$ at $t=\ell$.  Thus either  $v_{-j-1}v_{-j}v_{-j+1}$ or $v_{-j-2}v_{-j-1}v_{-j}v_{-j+1}$ is a path in $G$, and we can extend our trail to $v_{-j-1}v_{-j}\dots v_0$ or $v_{-j-2}v_{-j-1}v_{-j}\dots v_0$. 
 It should be obvious that no forcing edge will appear in this walk multiple times (by our construction). If $uv$ is not a forcing edge, then it can only appear in our walk if  $u'$ forced $u$ and $v$ forced $v'$. Let $u' u v v'$ be the first occurrence of $uv$ in our walk. 
 If $uv$ were to occur again, either $u$ or $v$ would need to be forced at this later time, but this cannot happen because $u$ and $v$ are both black at this point. \epf

\begin{ex}\label{extrail} {\rm Let $G$ be the graph in Figure \ref{longfish}.   As shown by the numbering in the figure, $\pt(G)=4$, but $G$ does not contain a path of length 4.  The trail produced by the method of proof used in  Lemma \ref{ztraillem} is $abcdecf$ and has length 6.}

  \begin{figure}[!ht] 
\begin{center}\scalebox{.4}{\includegraphics{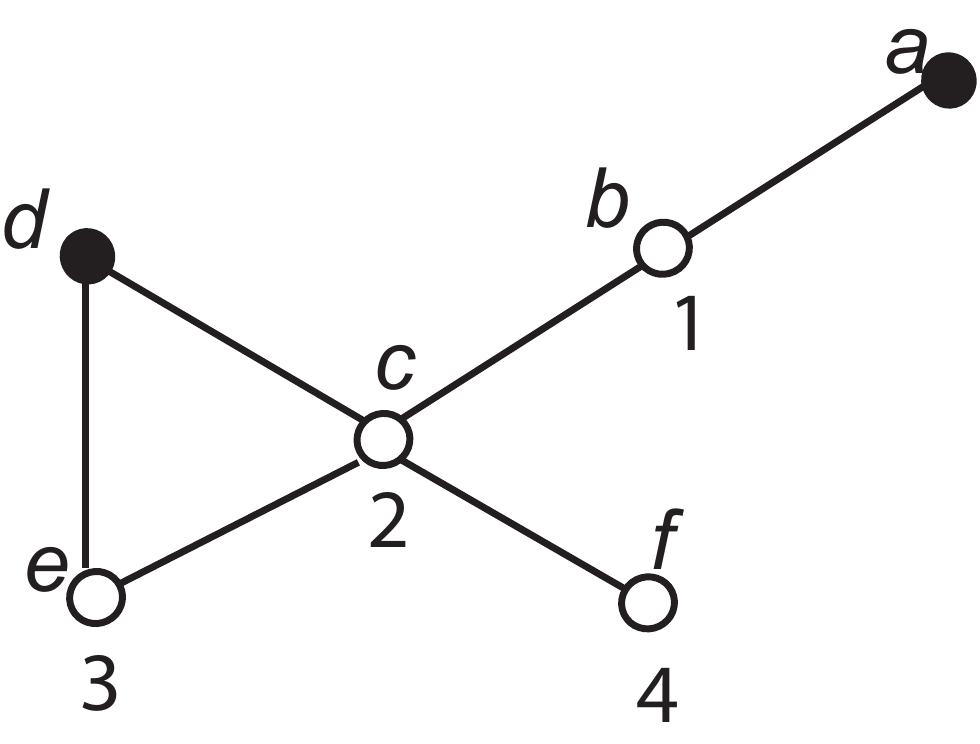}}
\caption{A graph $G$ that does not have a path of length $\pt(G)$} \label{longfish}
\end{center}
\end{figure}
\end{ex}

\begin{thm}\label{zdiamtree} Let $T$ be a tree and $B$ be a minimum zero forcing set of $T$. Then $\pt(T,B) \le \diam(T)$. Hence, $\pt(T) \le \PT(T) \le \diam(T)$.
\end{thm}

\bpf Choose $B$ to be a minimum zero forcing set such that $\pt(T,B) = \PT(T)$. By Lemma \ref{ztraillem}, there exists a trail in $T$ of length at least $\pt(T,B)$. Since  between any two vertices  in a tree  there is a unique path, any trail is a path and the diameter of $T$ must be the length of the longest path in $T$. Therefore, \[\pt(T) \le \PT(T) =\pt(T,B)  \le \diam(T).\qedhere\]
\epf

 The diameter of a graph $G$ can get arbitrarily larger than its minimum propagation time. The next example exhibits this result, but   first we observe that if
 $G$ is a graph having exactly $\ell$ leaves, then $\ZFN(G)\ge \lceil\frac \ell 2\rceil$ since at most two leaves can be on a maximal forcing chain.

\begin{figure}[!ht] 
\begin{center}\scalebox{.5}{\includegraphics{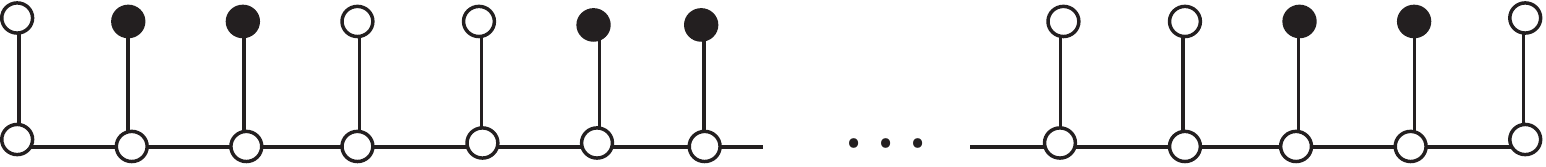}}
\caption{A $k$-comb} \label{comb}
\end{center}
\end{figure}

\begin{ex}{\rm 
 To construct a $k$-comb, we append a leaf to each vertex of a path on $k$ vertices, as shown in Figure \ref {comb} (our $k$-comb is the special case $P_{k,2}$ of a more general type of a comb $P_{k,\ell}$ defined in \cite{Sev}).   Let $G$ denote a $k$-comb where $k\equiv 0\mod 4$. It is clear that $\diam(G) = k+1$. 
 If we number the leaves in path order starting with one, then the set  $B$ consisting of every leaf whose number is congruent to 2 or 3$\mod 4$ (shown in black in the Figure $\ref{comb}$) is a zero forcing set, and  $\ord B = \frac{k}{2}$.  Since $\ZFN(G)\ge \frac k 2$, $B$ is a minimum zero forcing set.   Then $\pt(G) \leq \pt(G, B)= 3$. Since $|G|=2k$, $\ZFN(G)=\frac{k}{2}$, and $\pt(G)\ge \frac  {\ord G-\ZFN(G)} {\ZFN(G)}$, $\pt(G) \geq 3$. Therefore, $\pt(G) = 3$.  Thus the $\diam(G)=k+1$ is arbitrarily larger than $\pt(G)=3$. 
}\end{ex}


\noi{\bf Acknowledgement}
The authors thank the referees for many helpful suggestions.


\end{document}